\documentclass[12pt]{article}                
\usepackage{graphicx}
\usepackage{psfrag}
\usepackage{amsmath, amssymb}
\usepackage{mathptmx}
\usepackage{ntheorem}

\newcommand{\Av}{{\rm Av}}

\newcommand{\M}{{\cal M}}
\newcommand{\FF}{{\cal F}}

\newcommand{\II}{{\cal I}}

\newcommand{\Ind}{{\hspace{0.3mm}{\rm I}\hspace{0.1mm}}}
\newcommand{\eps}{{\varepsilon}}
\newcommand{\ch}{{\mbox{\rm ch}}}
\newcommand{\myth}{{\mbox{\rm th}}}
\newcommand{\smsp}{\hspace{0.3mm}}
\newcommand{\e}{\mathbb{E}}

\newcommand{\p}{\mathbb{P}}

\newcommand{\Reals}{\mathbb{R}}
\newcommand{\Natural}{\mathbb{N}}

\newcommand\qed{\hfill\hbox{\rlap{$\sqcap$}$\sqcup$}}

\newcommand{\sbar}{{\bar{s}}}
\newcommand{\txi}{{\tilde{\xi}}}

\newtheorem{lemma}{Lemma}
\newtheorem{theorem}{Theorem}

\theoremstyle{nonumberplain}
\newcommand\specialref{}

\begin{document}

\title{Structure of $1$-RSB asymptotic Gibbs measures\\ in the diluted $p$-spin models}
\author{Dmitry Panchenko\thanks{Dept. of Mathematics, Texas A\&M University, panchenk@math.tamu.edu. Partially supported by NSF grant.}\\
}
\date{}
\maketitle
\begin{abstract}
In this paper we study asymptotic Gibbs measures in the diluted $p$-spin models in the so called $1$-RSB case, when the overlap takes two values $q_*, q^*\in [0,1].$ When the external field is not present and the overlap is not equal to zero, we prove that such asymptotic Gibbs measures are described by the M\'ezard-Parisi ansatz conjectured in \cite{Mezard}. When the external field is present, we prove that the overlap can not be equal to zero and all $1$-RSB asymptotic Gibbs measures are described by the M\'ezard-Parisi ansatz. Finally, we give a characterization of the exceptional case when there is no external field and the smallest overlap value $q_*=0$, although it does not go as far as the M\'ezard-Parisi ansatz.  Our approach is based on the cavity computations combined with the hierarchical exchangeability of pure states.
\end{abstract} 
\vspace{0.5cm}

\section{Introduction and main result}

In \cite{Mezard}, M\'ezard and Parisi studied the diluted $p$-spin model for $p=2$ and described the structure of the Gibbs measure in the infinite-volume limit together with the corresponding formula for the free energy. They only formulated the $1$-step replica symmetry breaking ($1$-RSB) solution, but their ansatz has a natural extension to the general $r$-RSB case. It is expected that the same solution is valid for other diluted models as well, for example, for the random $K$-sat model and, possibly, for most mean field spin glass models. The origin of the M\'ezard-Parisi ansatz was partially explained in \cite{HEPS} via the hierarchical exchangeability of pure states combined with the hierarchical version of the Aldous-Hoover representation proved in \cite{AP}. However, as was also explained in \cite{HEPS}, some obstacles still remain in the form of additional symmetries between pure states, expressed by saying that `multi-overlaps between pure states are determined by their overlaps'. In this paper, we will prove the $1$-RSB M\'ezard-Parisi ansatz for diluted $p$-spin models in the case when the external field is present or when the overlap is not equal to zero. We will also show that the overlap can not be equal to zero in the presence of the external field. In the case when there is no external field and the smallest overlap value is zero, our approach will give information only about `odd moments' and we will not be able to recover the M\'ezard-Parisi ansatz completely. 

Most of our approach is rather general and can be extended to the $r$-RSB case, as well as to other models, such as the random $K$-sat model. However, the last step in the argument uses the special form of the $p$-spin model in a rather ad hoc way, and improving upon this could lead to progress in the general $r$-RSB case and for other diluted models. To understand the motivation for what we do in this paper, one should at least read the introduction in \cite{HEPS}, even though we will repeat all necessary definitions. In \cite{HEPS}, we used the random $K$-sat model as an example to illustrate the general approach, but the same results hold for the diluted $p$-spin models practically verbatim. The only place in \cite{HEPS} where the specific form of the $K$-sat model was used was in Lemma 1, where the self-averaging of the free energy was proved, and one can easily check that the same proof works for the diluted $p$-spin model. 

Consider an integer $p\geq 2$, the connectivity parameter $\lambda>0,$ the inverse temperature $\beta>0$ and the external field $h\in \Reals$. Consider a random function
\begin{equation}
\theta(\sigma_1,\ldots,\sigma_p)=\beta g \sigma_1\cdots \sigma_p
\label{Deftheta}
\end{equation}
on $\{-1,+1\}^p$, where $g$ is a standard Gaussian random variable. Let $(\theta_k)_{k\geq 1}$ be a sequence of independent copies of the function $\theta$, defined in terms of independent copies $(g_k)_{k\geq 1}$ of $g$. Then, using this sequence, the Hamiltonian $H_N(\sigma)$ of the diluted $p$-spin model on the space of spin configurations $\Sigma_N = \{-1,+1\}^N$ is defined by
\begin{equation}
H_N(\sigma) =  \sum_{k\leq \pi(\lambda N)}
\theta_k(\sigma_{i_{1,k}},\ldots,\sigma_{i_{p,k}})
+h\sum_{1\leq i\leq N} \sigma_i,
\label{Ham2}
\end{equation}
where $\pi(\lambda N)$ is a Poisson random variable with the mean $\lambda N$ and the indices $(i_{j,k})_{j,k\geq 1}$ are i.i.d. uniform on $\{1,\ldots,N\}$. The quantity\begin{equation}
F_N = \frac{1}{N}\e \log \sum_{\sigma\in \Sigma_N} \exp H_N(\sigma)
\end{equation}
is called the free energy of the model, and the probability measure on $\Sigma_N$ defined by 
\begin{equation}
G_N(\sigma) = \frac{1}{Z_N} \exp H_N(\sigma)
\label{GibbsN}
\end{equation} 
is called the Gibbs measure, where the normalizing factor $Z_N$ is called the partition function. The main goal in this model, as in other spin glass models, is to compute the limit of the free energy $F_N$ in the infinite-volume limit $N\to\infty$. In particular, any small perturbations of the Hamiltonian that do not affect the limit of the free energy are allowed, as long as they yield some useful information about the Gibbs measure. In this paper, we will utilize perturbations of two kinds to ensure that in the infinite-volume limit the Gibbs measure satisfies the Ghirlanda-Guerra identities and cavity equations. These perturbations will be reviewed in Section \ref{Sec2label}.

Before we state our main result, let us first recall the definition of asymptotic Gibbs measures introduced in \cite{Pspins} and also used in \cite{HEPS} (see \cite{Austin2} for a different approach via exchangeable random measures).

\medskip
\noindent
\textbf{Asymptotic Gibbs measures.}
Let $(\sigma^\ell)_{\ell\geq 1}$ be an i.i.d. sequence of replicas from the Gibbs measure $G_N$ and let $\mu_N$ be the joint distribution of the array of all spins on all replicas $(\sigma_i^\ell)_{1\leq i\leq N, \ell\geq 1}$ under the average product Gibbs measure $\e G_N^{\otimes \infty}$,
\begin{equation}
\mu_N\Bigl( \bigl\{\sigma_i^\ell = a_i^\ell \ :\ 1\leq i\leq N, 1\leq \ell \leq n \bigr\} \Bigr)
=
\e G_N^{\otimes n}\Bigl( \bigl\{\sigma_i^\ell = a_i^\ell \ :\ 1\leq i\leq N, 1\leq \ell \leq n \bigr\} \Bigr)
\label{muN}
\end{equation}
for any $n\geq 1$ and any $a_i^\ell \in\{-1,+1\}$.  We extend $\mu_N$ to a distribution on $\{-1,+1\}^{\Natural\times\Natural}$ simply by setting $\sigma_i^\ell=1$ for $i\geq N+1.$ Let $\M$ denote the set of all possible limits of $(\mu_N)$ over subsequences with respect to the weak convergence of measures on the compact product space $\{-1,+1\}^{\Natural\times\Natural}$. 

Notice that the distribution of the Hamiltonian (\ref{Ham2}) is invariant under the permutations of the coordinates of $\sigma$. Because of this property, called the symmetry between sites, all measures in $\M$ inherit from $\mu_N$ the invariance under the permutation of both spin and replica indices $i$ and $\ell.$ By the Aldous-Hoover representation \cite{Aldous}, \cite{Hoover2} for such distributions, for any $\mu\in\M$, there exists a measurable function $s:[0,1]^4\to\{-1,+1\}$ such that $\mu$ is the distribution of the array 
\begin{equation}
s_i^\ell=s(w,u_\ell,v_i,x_{i,\ell}),
\label{sigma}
\end{equation}
where the random variables $w,(u_\ell), (v_i), (x_{i,\ell})$ are i.i.d. uniform on $[0,1]$. The function $s$ is defined uniquely for a given $\mu\in \M$ up to measure-preserving transformations (Theorem 2.1 in \cite{Kallenberg}), so we can identify the distribution $\mu$ of array $(s_i^\ell)$ with $s$. Since $s$ takes values in $\{-1,+1\}$, the distribution $\mu$ can be encoded by the function
\begin{equation}
\sigma(w,u,v) = \e_x\, s(w,u,v,x),
\label{fop}
\end{equation}
where $\e_x$ is the expectation in $x$ only. The last coordinate $x_{i,\ell}$ in (\ref{sigma}) is independent for all pairs $(i,\ell)$, so it plays the role of `flipping a coin' with the expected value $\sigma(w,u_\ell,v_i)$. Therefore, given the function (\ref{fop}), we can redefine $s$ by
\begin{equation}
s(w,u_\ell,v_i,x_{i,\ell}) = 2 \Ind \Bigl(x_{i,\ell} \leq \frac{1+ \sigma(w,u_\ell,v_i) }{2}\Bigr) -1
\label{sigmatos}
\end{equation}
without affecting the distribution of the array $(s_i^\ell)$. 

We can also view the function $\sigma$ in (\ref{fop}) in a more geometric way as a random  measure on the space of functions, as follows. Let $du$ and $dv$ denote the Lebesgue measure on $[0,1]$ and let us define a (random) probability measure 
\begin{equation}
G = G_w = du \circ \bigl(u\to \sigma(w,u,\cdot)\bigr)^{-1}
\label{Gibbsw}
\end{equation}
on the space of functions of $v\in [0,1]$,
\begin{equation}
H = L^2\bigl([0,1], dv \bigr) \cap \bigl\{ \|\sigma\|_\infty \leq 1 \bigr\},
\label{spaceH}
\end{equation} 
equipped with the topology of $L^2([0,1], dv)$. We will denote by $\sigma^1\cdot \sigma^2$ the scalar product in $L^2([0,1], dv)$ and by $\|\sigma\|$ the corresponding $L^2$ norm. The random measure $G$ in (\ref{Gibbsw}) is called an {asymptotic Gibbs measure}. The whole process of generating spins can be broken into several steps:
\begin{enumerate}
\item[(i)] generate the Gibbs measure $G=G_w$ using the uniform random variable $w$; 

\item[(ii)] consider an i.i.d. sequence $\sigma^\ell = \sigma(w,u_{\ell},\cdot)$ of replicas from $G$, which are functions in $H$; 

\item[(iii)] plug in i.i.d. uniform random variables $(v_i)_{i\geq 1}$ to obtain the array $\sigma^\ell(v_i) = \sigma(w,u_\ell,v_i)$;

\item[(iv)] finally, use this array to generate spins as in (\ref{sigmatos}). 
\end{enumerate}
The M\'ezard-Parisi ansatz in \cite{Mezard} predicts that all asymptotic Gibbs measures (possibly, under a small perturbation of the Hamiltonian) have a very special structure. We will not repeat here what this structure is expected to be in general (see \cite{HEPS} for details) and will only describe it in the so called $1$-RSB case considered in \cite{Mezard}.

\medskip
\noindent
\textbf{The $1$-RSB M\'ezard-Parisi ansatz.} Suppose that an asymptotic Gibbs measure $G$ is such that, with probability one over the choice of this random measure, the scalar product $\sigma^1\cdot \sigma^2$ (also called the overlap) of points $\sigma^1$ and $\sigma^2$ in the support of $G$ can take one of the two non-random values $q_*< q^*$. In fact, this just means that the self-overlap is always $\sigma^1\cdot \sigma^1 = q^*$, so that the measure $G$ is supported on the sphere $\|\sigma\|^2 = q^*$, and the overlap of two different points is $\sigma^1\cdot \sigma^2=q_*$.  Of course, this also means that the measure $G$ is purely atomic,
\begin{equation}
G(\sigma_\alpha) = V_\alpha \ \mbox{ for }\ \alpha\in\Natural,
\label{Gibbsw1RSB}
\end{equation}
and we will assume that the atoms, which are called the pure states, are always enumerated in the decreasing order of their weights, $V_1 > V_2> \ldots > V_\alpha >\ldots$. For simplicity of notation, we will keep the dependence of the function $\sigma_\alpha$ and the weights $V_\alpha$ on $w$ implicit. Notice that in order to describe the distributions of all spins generated in steps (i) -- (iv) above, in the $1$-RSB case we need to describe the joint distribution of the weights $(V_\alpha)_{\alpha\in\Natural}$ and the array $(\sigma_\alpha(v_i))_{\alpha, i\in\Natural}$. The $1$-RSB M\'ezard-Parisi ansatz predicts the following.
\begin{enumerate}
\item[(a)] The weights $(V_\alpha)_{\alpha\in\Natural}$ and the array $(\sigma_\alpha(v_i))_{\alpha, i\in\Natural}$ are independent. 

\item[(b)] The weights $(V_\alpha)_{\alpha\in\Natural}$ have the Poisson-Dirichlet distribution $PD(\zeta)$ for some $\zeta\in (0,1)$.

\item[(c)] There exists a function $f:[0,1]^3 \to [-1,1]$ such that
\begin{equation}
\bigl(\sigma_\alpha(v_i) \bigr)_{\alpha, i\in\Natural}
\, \stackrel{d}{=}\,
\bigl( f(\omega, \omega^i, \omega_\alpha^i) \bigr)_{\alpha, i\in\Natural},
\label{sigmaf}
\end{equation}
where all $\omega, \omega^i, \omega_\alpha^i$ are i.i.d. random variables with the uniform distribution on $[0,1]$. 
\end{enumerate}
Let us discuss these properties in more detail. First of all, when we sample replicas $(\sigma^\ell)$ from the Gibbs measure $G$, we sample them from the list of pure states $(\sigma_\alpha)_{\alpha\in \Natural}$ in $H$ according to weights $(V_\alpha)_{\alpha\in\Natural}$, which have the Poisson-Dirichlet distribution $PD(\zeta)$. We remind that if $(x_\alpha)_{\alpha \in \Natural}$ is the decreasing enumeration of a Poisson process on $(0,\infty)$ with the mean measure $\zeta x^{-1-\zeta}dx$ for some $\zeta\in (0,1)$ then the distribution of the sequence
\begin{equation}
V_\alpha = \frac{x_\alpha}{\sum_{\alpha \geq 1} x_\alpha} 
\end{equation}
is called the {Poisson-Dirichlet distribution} $PD(\zeta)$. It is well known (see e.g. Section 2.2 in \cite{SKmodel}) that the parameter 
$$
\zeta = 1-  \e \sum_{\alpha\geq 1} V_\alpha^2 = \e \sum_{\alpha\not =\beta} V_\alpha V_{\beta}
$$
represents the probability that two pure states sampled according to $(V_\alpha)$ will be different. Then, independently from the weights of the pure states, we generate the array $(\sigma_\alpha(v_i))_{\alpha, i\in \Natural}$ as in (\ref{sigmaf}). The random variable $\sigma_\alpha(v_i)$ is called the magnetization of the $i^{\mathrm{th}}$ spin inside the pure state $\alpha$ and, for a fixed $\omega$, the function $f(\omega, \,\cdot\,,\,\cdot\,)$ in (\ref{sigmaf}) represents the functional order parameter of the M\'ezard-Parisi ansatz. Conditionally on this auxiliary randomness $\omega$, the spin magnetizations are generated independently over $i\geq 1$ and, for each $i$, are generated in a completely symmetric exchangeable fashion over the pure states $\alpha\geq 1$. For example, (\ref{sigmaf}) implies that the multi-overlaps
$$
\int\! \sigma_{\alpha_1}(v)\cdots \sigma_{\alpha_n}(v)  \, dv
$$ 
of pure states $\alpha_1,\ldots,\alpha_n$ (not necessarily all different) are equal in distribution to
$$
\e_i f(\omega, \omega^i, \omega_{\alpha_1}^i) \cdots f(\omega, \omega^i, \omega_{\alpha_n}^i),
$$
where $\e_i$ denotes the expectation in the random variables $\omega^i, (\omega_{\alpha}^i)_{\alpha\geq 1}$. Obviously, this quantity depends only on the values
$$
\Ind(\alpha_\ell = \alpha_{\ell'}) = \Ind(\sigma_{\alpha_\ell} \cdot \sigma_{\alpha_{\ell'}} = q^*)
\ \mbox{ for }\
1\leq \ell,\ell'\leq n
$$
determined by the overlaps between pure states, so, in other words, multi-overlaps are determined by the overlaps.

We will prove the $1$-RSB M\'ezard-Parisi ansatz under a small perturbation of the Hamiltonian (\ref{Ham2}). In the next section, we will define a slightly modified Hamiltonian 
\begin{equation}
H_N^{\mathrm{pert}}(\sigma) = H_N(\sigma)+ h_N^{\mathrm{pert}}(\sigma)
\label{HNpert}
\end{equation}
for some small perturbation $h_N^{\mathrm{pert}}(\sigma)$ that does not affect the limit of the free energy and, from now on, consider asymptotic Gibbs measures corresponding to this perturbed Hamiltonian. The perturbation will force the asymptotic Gibbs measures to satisfy several properties sufficient to prove the following.
\begin{theorem}\label{Th1}
If $h=0$ then any $1$-RSB asymptotic Gibbs measure such that $q_* \not = 0$ satisfies the M\'ezard-Parisi ansatz. If $h\not =0$ then $q_*\not = 0$ and any $1$-RSB asymptotic Gibbs measure satisfies the M\'ezard-Parisi ansatz. 
\end{theorem} 
In the next section, we will also complement Theorem \ref{Th1} and explain what happens when $h=0$ and $q_*=0$. We should also mention that, in general, Theorem \ref{Th1} by itself does not say anything about the free energy. However, if $h\not = 0$ and one could show that, in some region of parameters $(\lambda,\beta)$, all asymptotic Gibbs measures are $1$-RSB then one could also recover the M\'ezard-Parisi $1$-RSB formula for the free energy when $p\geq 2$ is even, using \cite{FL, PT}. For the case $h=0$, it would be sufficient to show that for all small enough $h\not = 0$, all asymptotic Gibbs measures are $1$-RSB. Then one could also recover the formula for the free energy by letting $h$ go to zero. 

In the next section, we will describe two kinds of perturbation of the Hamiltonian and the corresponding properties they ensure---some consequences of the Ghirlanda-Guerra identities and the cavity equations. In Section \ref{Sec3label}, we will rewrite the cavity equations specifically for the $1$-RSB case and in Section \ref{Sec4label}, using the properties of the Poisson-Dirichlet distribution of the pure state weights, we will deduce a variant of the cavity equations for the pure states. In Section \ref{Sec5label}, we will  prove the key consequence of the cavity equations, and in Section \ref{Sec6label} we will use it to prove Theorem \ref{Th1} in the case when $h=0$. In Section \ref{Sec7label}, we will study the case when $q_*=0$, and in Section \ref{Sec8label} we will prove Theorem \ref{Th1} in the case when $h\not=0$.

\section{Properties of Gibbs measures via perturbations}\label{Sec2label}

The perturbation term $h_N^{\mathrm{pert}}(\sigma)$ in (\ref{HNpert}) will consist of two parts,
\begin{equation}
h_N^{\mathrm{pert}}(\sigma)= h_N^{\mathrm{1}}(\sigma)+ h_N^{\mathrm{2}}(\sigma).
\label{hNparts}
\end{equation}
Each part will be responsible for a certain property of the asymptotic Gibbs measures. 

\medskip
\noindent
\textbf{Perturbation of the first kind.}
For each $\ell\geq 1$, let us consider the process $g_{N,\ell}(\sigma)$ on  $\{-1,+1\}^N$ given by
\begin{equation}
g_{N,\ell}(\sigma)
=
\frac{1}{N^{\ell/2}}
\sum_{1\leq i_1,\ldots,i_\ell \leq N} g_{i_1,\ldots, i_\ell} \sigma_{i_1}\ldots\sigma_{i_\ell},
\label{mixedppert}
\end{equation}
where $(g_{i_1,\ldots, i_\ell})$ are i.i.d. standard Gaussian random variables, and define
\begin{equation}
h_N^{\mathrm{1}}(\sigma) = s_N \sum_{\ell\geq 1} 2^{-\ell} x_\ell^N \smsp g_{N,\ell}(\sigma),
\label{mixedHpert}
\end{equation}
where $s_N = N^{\gamma}$ for any $\gamma\in (1/4, 1/2)$ and parameters $x_\ell^N\in[0,3]$ for all $\ell\geq 1$. In Section 2 in \cite{HEPS} it was explained that this perturbation does not affect the limit of the free energy and, for some choice of parameters $x^N = (x_\ell^N)_{\ell\geq 1}$, all asymptotic Gibbs measures satisfy the Ghirlanda-Guerra identities \cite{GG}. We will not repeat the definition of the Ghirlanda-Guerra identities here and will only mention their main consequences proved in Theorem 1 in \cite{HEPS} (more precisely, the consequences of the invariance principle discovered in \cite{PUltra} that follows from the Ghirlanda-Guerra identities). Namely, Theorem 1 and the discussion right after the Corollary 1 in \cite{HEPS} imply that any $1$-RSB asymptotic Gibbs measure satisfies the properties (a) and (b) in the $1$-RSB M\'ezard-Parisi ansatz and property (c) is replaced with
\begin{enumerate}
\item[(c)${}^\prime$] There exists a function $f:[0,1]^4 \to [-1,1]$ such that
\begin{equation}
\bigl(\sigma_\alpha(v_i) \bigr)_{\alpha, i\in\Natural}
\, \stackrel{d}{=}\,
\bigl( f(\omega, \omega_\alpha, \omega^i, \omega_\alpha^i) \bigr)_{\alpha, i\in\Natural},
\label{sigmaf2}
\end{equation}
where all $\omega, \omega_\alpha, \omega^i, \omega_\alpha^i$ are i.i.d. random variables with the uniform distribution on $[0,1]$. 
\end{enumerate}
This means that our main goal now is to show that we can  replace $f$ on the right hand side of (\ref{sigmaf2}) by a function that does not depend on $\omega_\alpha$, proving the representation (\ref{sigmaf}) that encodes a much simpler and much more symmetric structure than (\ref{sigmaf2}). As we already mentioned above, in the case when $q_*=0$, we will not be able to prove Theorem \ref{Th1} and, instead, give the following characterization.
\begin{theorem}\label{Th1comp}
For almost all $(\omega, \omega_\alpha, \omega^i)$, the conditional distribution of $ f(\omega, \omega_\alpha, \omega^i, \omega_\alpha^i)$ in (\ref{sigmaf2}) given $(\omega, \omega_\alpha, \omega^i)$ is symmetric if and only if $q_* = 0$.
\end{theorem} 
Both Theorem \ref{Th1} and Theorem \ref{Th1comp} will be deduced from (\ref{sigmaf2}) and the cavity equations that can be proved with the help of the following perturbation. 

\medskip
\noindent
\textbf{Perturbation of the second kind.} Consider a sequence $(c_N)$ such that $c_N \uparrow \infty$ and $|c_{N+1}-c_N|\to 0$.  Consider an i.i.d. sequence of indices $(i_{j,k,\ell})_{j,k,\ell\geq 1}$ with the uniform distribution on $\{1,\ldots,N\}$, let $\pi(c_N)$ be a Poisson random variable with the mean $c_N,$ $(\pi_\ell(\lambda p))_{\ell\geq 1}$ be i.i.d. Poisson with the mean $\lambda p$, and $(\theta_{\ell,k})_{\ell,k\geq 1}$ are i.i.d. copies of the function (\ref{Deftheta}). We define the second perturbation term by
\begin{equation} 
h_N^{\mathrm{2}}(\sigma)
=
\sum_{\ell \leq \pi(c_N)} 
\log \Av \exp \Bigl(\sum_{k\leq \pi_\ell(\lambda p)}
\theta_{\ell,k} (\sigma_{i_{1,k,\ell}}, \ldots, \sigma_{i_{p-1,k,\ell}}, \eps)
+h\eps\Bigr),
\label{Hampert}
\end{equation}
where $\Av$ denotes the average over $\eps\in\{-1,+1\}$. Notice that the condition $|c_{N+1}-c_N|\to 0$ implies that $c_N/N \to 0$ and, therefore, this perturbation also does not affect the limit of the free energy. The perturbation (\ref{Hampert}) was introduced in \cite{Pspins} in order to prove the cavity equations for the spins. In a few words, the main idea behind this perturbation is that it represents the affect on the Hamiltonian of adding $\pi(c_N)$ spins to the system and treating them as cavity coordinates. This adds some stability to the Gibbs measure when we consider a finite number of additional coordinates as cavity coordinates (they are lost inside the big crowd of $\pi(c_N)$ cavity coordinates, so to speak) and this stability allows one to prove the following cavity equations. 

We will need to pick various sets of different spin coordinates in the array $(s_i^\ell)$ in (\ref{sigma}), and it is inconvenient to enumerate them using one index $i\geq 1$. Instead, we will use multi-indices $I= (i_1,\ldots, i_n)$ for $n\geq 1$ and $i_1,\ldots, i_n\geq 1$ and consider 
\begin{equation}
s_{I}^\ell = s(w,u_\ell, v_{I},x_{I,\ell}),
\end{equation}
where all the coordinates are uniform on $[0,1]$ and independent over different sets of indices. For convenience, below we will separate averaging with respect to random variables that are indexed by different replica indices $\ell$, and for this purpose we will use the notation
\begin{equation}
s_{I} = s(w,u, v_{I},x_{I}).
\label{s}
\end{equation}
Now, take arbitrary integers $n, m, q\geq 1$ such that $n\leq m.$ The index $q$ will represent the number of replicas selected, $m$ will be the total number of spin coordinates and $n$ will be the number of cavity coordinates. For each replica index $\ell\leq q$ we  consider an arbitrary subset of coordinates $C_\ell\subseteq \{1,\ldots, m\}$ and split them into cavity and non-cavity coordinates,
\begin{equation}
C_\ell^1 = C_\ell\cap \{1,\ldots, n\},\,\,\,
C_\ell^2=C_\ell\cap \{n+1,\ldots,m\}.
\label{C12}
\end{equation}
The following quantities represent the cavity fields for $i\geq 1$,
\begin{equation}
A_i(\eps)=\sum_{k\leq \pi_i(\lambda p)} \theta_{i,k}( s_{1,i,k}, \ldots, s_{p-1,i,k}, \eps) + h\eps,
\label{Ai}
\end{equation}
where $\eps\in \{-1,+1\}$, $(\pi_i(\lambda p))_{i\geq 1}$ are i.i.d. Poisson random variables with the mean $\lambda p$, and $(\theta_{i,k})_{k,i\geq 1}$ are i.i.d. copies of the function (\ref{Deftheta}). Let $\e'$ denote the expectation in $u$ and the random variables $x_I$ for all multi-indices $I$, and $\Av$ denote the uniform average over $(\eps_i)_{i\geq 1}$ in $\{-1,+1\}^\Natural$. Define
\begin{equation}
U_\ell =\, \e' \Av \prod_{i\in C_\ell^1} \eps_i \prod_{i\in C_\ell^2} s_i
\,\exp \sum_{i\leq n} A_i(\eps_i) \ \mbox{ and } \
V =\, \e' \Av \exp \sum_{i\leq n} A_i(\eps_i).
\label{Ul}
\end{equation}
Then Theorem $1$ in \cite{Pspins} states that, for any asymptotic Gibbs measure, we have
\begin{equation}
\e \prod_{\ell\leq q} \e' \prod_{i\in C_\ell}s_i
=\e \prod_{\ell\leq q} \frac{U_\ell}{V}.
\label{SC}
\end{equation}
The left hand side can be written using replicas as $\e \prod_{\ell\leq q} \prod_{i\in C_\ell}s_i^\ell$, so it represent an arbitrary joint moment of spins in the array (\ref{sigma}). The right hand side expresses what happens to this joint moment in the infinite-volume limit when we treat the first $n$ spins as cavity coordinates. We will utilize these cavity equations to show that the function $f$ in (\ref{sigmaf2}) can be replaced by a function in (\ref{sigmaf}) that does not depend on the coordinate $\omega_{\alpha}$.
 
Let us remark that the proof of the cavity equations in \cite{Pspins} was given only in the case when $h=0$, but it is identical in the case when $h\not =0$. Simply, the cavity field has one additional term $h\eps$. Also, the perturbations (\ref{mixedHpert}) and (\ref{Hampert}) were considered in \cite{HEPS} and \cite{Pspins} separately and not at the same time as we do here. However, it is not difficult to see by inspecting the proofs there that these two perturbations do not interfere with each other and we can obtain all the corresponding consequences for the asymptotic Gibbs measures at the same time, i.e. (a), (b), $\mbox{(c)}^\prime$ and (\ref{SC}). For example, since the perturbation (\ref{mixedHpert}) is of a smaller order, its affect on the cavity fields will be negligible and can be ignored in the proof of the cavity equations (\ref{SC}). On the other hand, since the perturbation (\ref{Hampert}) is also of a smaller order, it does not affect the self-averaging of the free energy, which was the main reason behind the Ghirlanda-Guerra identities.

\section{Rewriting the cavity equations}\label{Sec3label}

In this section, we will rewrite the cavity equations (\ref{SC}) in the $1$-RSB case using the discrete nature of the Gibbs measure in (\ref{Gibbsw1RSB}) and the representation of spin magnetizations inside the pure states stated in (\ref{sigmaf2}). This is nothing but a straightforward reformulation in a couple of steps. In the first step, it will be convenient to extend the definition of the function $\theta$ in (\ref{Deftheta}) from $\{-1,+1\}^p$ to $[-1,1]^p$ as follows. Since the product $\sigma_1\cdots \sigma_p$ in (\ref{Deftheta}) takes only two values $\pm 1$, we can write
\begin{equation}
\exp \theta(\sigma_1,\ldots,\sigma_p)= 
\ch(\beta g)\bigl(1+ \myth(\beta g)\, \sigma_1\cdots \sigma_p \bigr).
\label{expTheta}
\end{equation}
In a moment, we will be averaging $\exp \theta$ over the coordinates $\sigma_1,\ldots,\sigma_p$ independently of each other, so the resulting average will be of the same form with $\sigma_i$ taking values in $[-1,1].$ We will again represent this average as $\exp \theta$ with $\theta$ now defined by
\begin{equation}
\theta(\sigma_1,\ldots,\sigma_p)=\log\Bigl(\ch(\beta g)\bigl(1+ \myth(\beta g)\, \sigma_1\cdots \sigma_p \bigr)\Bigr).
\label{DefthetaG}
\end{equation}
Of course, on the set $\{-1,+1\}^p$ this definition coincides with (\ref{Deftheta}). 

Let us write $\e' = \e_u \e_x$, where $\e_u$ denote the expectation in $u$ and $\e_x$ denotes the expectations in the random variables $x_I$ for all multi-indices $I$. Recalling (\ref{fop}) and (\ref{s}), we can write
\begin{equation}
\sbar_{I} := \e_x s_I = \sigma(w,u, v_{I}).
\label{sbar}
\end{equation}
If, similarly to (\ref{Ai}), we denote
\begin{equation}
\bar{A}_i(\eps)=\sum_{k\leq \pi_i(\lambda p)} \theta_{i,k}( \sbar_{1,i,k}, \ldots, \sbar_{p-1,i,k}, \eps) +h\eps
\label{Aibar}
\end{equation}
then, using (\ref{expTheta}) and (\ref{DefthetaG}), we can write
\begin{equation}
\e_x \exp \sum_{i\leq n}A_i(\eps)
=
\exp \sum_{i\leq n}\bar{A}_i(\eps).
\end{equation}
Therefore, if similarly to (\ref{Ul}) we define
\begin{equation}
\bar{U}_\ell =\, \e_u \Av \prod_{i\in C_\ell^1} \eps_i \prod_{i\in C_\ell^2} \sbar_i
\,\exp \sum_{i\leq n} \bar{A}_i(\eps_i) \ \mbox{ and } \
\bar{V} =\, \e_u \Av \exp \sum_{i\leq n} \bar{A}_i(\eps_i)
\label{Ulbar}
\end{equation}
then the cavity equations (\ref{SC}) can be rewritten as
\begin{equation}
\e \prod_{\ell\leq q} \e_u \prod_{i\in C_\ell} \sbar_i
=\e \prod_{\ell\leq q} \frac{\bar{U}_\ell}{\bar{V}}.
\label{SCbar}
\end{equation}
Simply, we averaged out the random variables $x_I$. Next, let us denote
\begin{equation}
\bar{A}_i = \log \Av \exp \bar{A}_i(\eps) \ \mbox{ and }\
\xi_i = \frac{\Av \eps \exp \bar{A}_i(\eps) }{\Av \exp \bar{A}_i(\eps)} 
= \frac{\Av \eps \exp \bar{A}_i(\eps) }{\exp \bar{A}_i}.
\end{equation}
Then, (\ref{Ulbar}) can be rewritten as
\begin{equation}
\bar{U}_\ell =\, \e_u \prod_{i\in C_\ell^1} \xi_i \prod_{i\in C_\ell^2} \sbar_i
\,\exp \sum_{i\leq n} \bar{A}_i \ \mbox{ and } \
\bar{V} =\, \e_u \exp \sum_{i\leq n} \bar{A}_i.
\label{Ulbar2}
\end{equation}
Finally, comparing the definition of the measure $G$ in (\ref{Gibbsw}) with the fact that in the $1$-RSB case the measure $G$ is discrete as in (\ref{Gibbsw1RSB}), the expectation $\e_u$ in $u$ corresponds to averaging over the points $\sigma_\alpha$ in the support of $G$ with the weights $V_\alpha$. We will use this observation simultaneously with the property (\ref{sigmaf2}). Therefore, if we now define
\begin{align}
s_I^\alpha &\ =\ f(\omega, \omega_\alpha, \omega^I, \omega_\alpha^I),
\label{sialpha}
\\
A_i^\alpha(\eps)&\ =\ \sum_{k\leq \pi_i(\lambda p)} \theta_{i,k}( s_{1,i,k}^\alpha, \ldots, s_{p-1,i,k}^\alpha, \eps) +h\eps,
\label{Aieps}
\\
A_i^\alpha &\ =\ \log \Av \exp A_i^\alpha(\eps),
\label{Aialpha}
\\
\xi_i^\alpha &\ =\ \frac{\Av \eps \exp A_i^\alpha(\eps) }{\Av \exp A_i^\alpha(\eps)}
= \frac{\Av \eps \exp A_i^\alpha(\eps) }{\exp A_i^\alpha},
\label{Sec3xiialpha}
\end{align}
and let $A^\alpha = \sum_{i\leq n} A_i^\alpha$ then (\ref{Ulbar}) can be redefined by (using equality in distribution (\ref{sigmaf2}))
\begin{equation}
\bar{U}_\ell =\, \sum_{\alpha\geq 1} V_\alpha \prod_{i\in C_\ell^1} \xi_i^\alpha \prod_{i\in C_\ell^2} s_i^\alpha\,\exp A^\alpha \ \mbox{ and } \
\bar{V} =\, \sum_{\alpha\geq 1} V_\alpha \exp A^\alpha.
\label{Ulbaralpha}
\end{equation}
Moreover, if we denote
\begin{equation}
V_\alpha' = \frac{V_\alpha \exp A^\alpha}{\bar{V}}
= \frac{V_\alpha \exp A^\alpha}{\sum_{\alpha\geq 1} V_\alpha \exp A^\alpha}
\label{Valpha}
\end{equation}
then the cavity equations (\ref{SCbar}) take form
\begin{equation}
\e \prod_{\ell\leq q} \sum_{\alpha\geq 1} V_\alpha \prod_{i\in C_\ell} s_i^\alpha
=\e \prod_{\ell\leq q} \sum_{\alpha\geq 1} V_\alpha' \prod_{i\in C_\ell^1} \xi_i^\alpha \prod_{i\in C_\ell^2} s_i^\alpha.
\label{SCbaralpha}
\end{equation}
We can also write this  as
\begin{equation}
\e \sum_{\alpha_1,\ldots, \alpha_q} V_{\alpha_1}\cdots V_{\alpha_q} \prod_{\ell\leq q}
\prod_{i\in C_\ell} s_i^{\alpha_\ell}
=\e \sum_{\alpha_1,\ldots, \alpha_q} V_{\alpha_1}' \cdots V_{\alpha_q}' \prod_{\ell\leq q} \prod_{i\in C_\ell^1} \xi_i^{\alpha_\ell} \prod_{i\in C_\ell^2} s_i^{\alpha_\ell}.
\label{SCnew}
\end{equation}
In the next section, we will use this form of the cavity equations to obtain a different form directly for the pure states that does not involve averaging over the pure states.

\section{Cavity equations for the pure states}\label{Sec4label}

Let $\FF$ be the $\sigma$-algebra generated by the random variables $g_{i,k}$ (the Gaussian coefficients of the functions $\theta_{i,k}$), $\pi_i(\lambda p), \omega, \omega^I$ for various indices, excluding only the random variables $\omega_{\alpha}$ and $\omega_{\alpha}^I$ that are indexed by $\alpha$. Conditionally on $\FF$, let $(\txi_i^{\alpha})_{i\leq n}$ be random vectors independent over $\alpha\geq 1$ with the distribution of $(\xi_i^{\alpha})_{i\leq n}$ in (\ref{Sec3xiialpha}) under the change of density
\begin{equation}
R^\alpha : = \frac{\exp \zeta A^\alpha}{ \e_\alpha \exp \zeta A^\alpha},
\label{Ralpha}
\end{equation}
where $\e_\alpha$ denotes the expectation in the random variables $\omega_{\alpha}$ and $\omega_{\alpha}^I$. Notice that this distribution does not depend on $\alpha$ so, conditionally of $\FF$, $(\txi_i^{\alpha})_{i\leq n}$ are i.i.d. for $\alpha\geq 1$. We will prove the following.
\begin{theorem}\label{Sec4Th}
The equality in distribution holds (not conditionally on $\FF$),
\begin{equation}
\bigl(\txi_i^{\alpha} \bigr)_{i\leq n,\alpha\in \Natural}
\stackrel{d}{=}
\bigl(s_i^{\alpha} \bigr)_{i\leq n,\alpha\in \Natural}.
\label{Sec4ThEq}
\end{equation}
\end{theorem}
\textbf{Proof.} 
We begin by noticing that the property (\ref{sigmaf2}) implies that the overlap of two pure states
\begin{equation}
R_{\alpha, \beta}:= \sigma_\alpha\cdot \sigma_\beta 
= \int_0^1\! \sigma_\alpha(v) \sigma_\beta(v) \, dv
\stackrel{d}{=} \e_i f(\omega, \omega_\alpha, \omega^i, \omega_\alpha^i) f(\omega, \omega_\beta, \omega^i, \omega_\beta^i),
\label{Sec3eq1}
\end{equation}
where $\e_i$ denotes the expectation in the random variables that depend on the spin index $i$. By the $1$-RSB assumption, $R_{\alpha, \beta}=q_*$ for $\alpha\not =\beta$ and $R_{\alpha, \beta}=q^*$ for $\alpha=\beta$. If we recall the definition (\ref{sialpha}), this implies that
\begin{equation}
R_{\alpha,\beta} =
\e_i \, s_i^{\alpha} s_i^{\beta}=q_*\Ind(\alpha\not =\beta) 
+ q^*\Ind(\alpha=\beta).
\label{Rab}
\end{equation}
In the cavity equations (\ref{SCnew}), let us now make a special choice of the sets $C_\ell^2$. For each pair $(\ell,\ell')$ of replica indices such that $1\leq \ell<\ell'\leq q$, take any integer $n_{\ell,\ell'}\geq 0$ and consider a set $C_{\ell,\ell'}\subseteq \{n+1,\ldots,m\}$ of cardinality $|C_{\ell,\ell'}|=n_{\ell,\ell'}$. Let all these sets be disjoint, which can be achieved by taking $m=n+\sum_{1\leq \ell<\ell'\leq q} n_{\ell,\ell'}.$ For each $\ell\leq q$, let
$$
C_\ell^2 = \Bigl(\bigcup_{\ell'>\ell} C_{\ell,\ell'}\Bigr) \bigcup \Bigl(\bigcup_{\ell'<\ell} C_{\ell',\ell}\Bigr).
$$
Then a given spin index $i\in \{n+1,\ldots,m\}$ appears in exactly two sets, say, $C_\ell^2$ and $C_{\ell'}^2$, and the expectation of (\ref{SCnew}) in $\omega^i,\omega_{\alpha_\ell}^i,\omega_{\alpha_{\ell'}}^i$ will produce a factor $\e_i \,s_i^{\alpha_\ell} s_i^{\alpha_{\ell'}} = R_{\alpha_\ell,\alpha_{\ell'}}$. For each pair $(\ell,\ell')$, there will be exactly $n_{\ell,\ell'}$ such factors, so averaging in (\ref{SCnew}) in the random variables $\omega^i,\omega_{\alpha_\ell}^i,\omega_{\alpha_{\ell'}}^i$ for all $i\in \{n+1,\ldots,m\}$ will result in
\begin{equation}
\e \sum_{\alpha_1,\ldots, \alpha_q} V_{\alpha_1}\cdots V_{\alpha_q} \prod_{\ell<\ell'} R_{\alpha_\ell, \alpha_{\ell'}}^{n_{\ell,\ell'}}
\prod_{\ell\leq q} \prod_{i\in C_\ell^1} s_i^{\alpha_\ell}
=\e \sum_{\alpha_1,\ldots, \alpha_q} V_{\alpha_1}' \cdots V_{\alpha_q}'  \prod_{\ell<\ell'} R_{\alpha_\ell, \alpha_{\ell'}}^{n_{\ell,\ell'}} \prod_{\ell\leq q} \prod_{i\in C_\ell^1} \xi_i^{\alpha_\ell}.
\label{SCagain}
\end{equation}
Approximating by polynomials, we can replace $\prod_{\ell<\ell'} R_{\alpha_\ell, \alpha_{\ell'}}^{n_{\ell,\ell'}}$ by an indicator of the set
\begin{equation}
C = \bigl\{(\alpha_1,\ldots, \alpha_q) \ | \ R_{\alpha_\ell, \alpha_{\ell'}} = q_{\ell,\ell'} \mbox{ for all } 1\leq \ell<\ell' \leq q\bigr\}
\end{equation}
for any choice of constraints $q_{\ell,\ell'}$ taking values $q_*$ or $q^*.$ In fact, since the overlaps take only two values, we can write this indicator as a finite linear combination of monomials with $n_{\ell,\ell'}$ taking values $0$ or $1$. We can also write the set $C$ as
\begin{equation}
C = \bigl\{(\alpha_1,\ldots, \alpha_q) \ |\ \alpha_\ell = \alpha_{\ell'} \mbox{ if and only if } q_{\ell,\ell'} =q^*\bigr\}.
\end{equation}
Therefore, (\ref{SCagain}) implies
\begin{equation}
\sum_{(\alpha_1,\ldots, \alpha_q)\in C} \e V_{\alpha_1}\cdots V_{\alpha_q} \prod_{\ell\leq q} \prod_{i\in C_\ell^1} s_i^{\alpha_\ell}
= \sum_{(\alpha_1,\ldots, \alpha_q)\in C} \e V_{\alpha_1}' \cdots V_{\alpha_q}'  \prod_{\ell\leq q} \prod_{i\in C_\ell^1} \xi_i^{\alpha_\ell}.
\label{SCF}
\end{equation}
Using the property (a) of the M\`erard-Parisi ansatz, which as we mentioned is the consequence of the perturbation of the first kind, we can rewrite the left hand side as
$$
\sum_{(\alpha_1,\ldots, \alpha_q)\in C} \e V_{\alpha_1}\cdots V_{\alpha_q} \,
\e \prod_{\ell\leq q} \prod_{i\in C_\ell^1} s_i^{\alpha_\ell}.
$$
Moreover, it is obvious from the definition of the array $s_i^\alpha$ in (\ref{sialpha}) that the second expectation does not depend on $(\alpha_1,\ldots, \alpha_q)\in C$. 

On the other hand, on the right hand side of (\ref{SCF}) both $V_\alpha'$ and $\xi_i^\alpha$ depend on the same random variables through the function $A_i^\alpha(\eps)$. However, the fact that by the property (b) of the M\'ezard-Parisi ansatz the sequence of weights $(V_\alpha)$ has the Poisson-Dirichlet distribution $PD(\zeta)$ allows us to overcome this obstacle as follows. We will consider the random variables $V_\alpha'$ and $\xi_i^\alpha$ as functions of $\omega_{\alpha}$ and $\omega_{\alpha}^I$ for various multi-indices $I$, conditionally on the $\sigma$-algebra $\FF$ defined above the equation (\ref{Ralpha}). Notice that, conditionally of $\FF$, the random pairs $(A^\alpha, (\xi_i^\alpha)_{i\leq n})$ are i.i.d. over $\alpha\geq 1$. Let $\rho:\Natural \to \Natural$ be the map that rearranges the weights $V_\alpha'$ in (\ref{Valpha}) in the decreasing order,
$$
V_{\rho(1)}' > V_{\rho(2)}' > \ldots > V_{\rho(\alpha)}' > \ldots.
$$
Then Theorem 2.6 in \cite{SKmodel} implies that
\begin{equation}
\Bigl(V_{\rho(\alpha)}' , \bigl(\xi_i^{\rho(\alpha)} \bigr)_{i\leq n} \Bigr)_{\alpha\geq 1}
\stackrel{d}{=}
\Bigl(V_{\alpha} , \bigl(\txi_i^{\alpha} \bigr)_{i\leq n} \Bigr)_{\alpha\geq 1},
\label{Sec4import}
\end{equation}
where the two sequences on the right hand side over $\alpha\geq 1$ are independent, the sequence $(V_\alpha)_{\alpha\geq 1}$ has the Poisson-Dirichlet distribution $PD(\zeta)$, the random vectors $(\txi_i^{\alpha})_{i\leq n}$ are i.i.d. over $\alpha\geq 1$ and have the distribution of $(\xi_i^{\alpha})_{i\leq n}$ under the change of density (\ref{Ralpha}). Since the distribution of the weights, $PD(\zeta)$, does not depend on the condition, the two sequences are also independent unconditionally. Together with the fact that $\rho$ is a bijection (this is a consequence of Theorem 2.6 in \cite{SKmodel} and is explained below equation (2.24) in \cite{SKmodel}) and $(\rho(\alpha_1),\ldots, \rho(\alpha_q))\in C$ if and only if $(\alpha_1,\ldots, \alpha_q)\in C$, the equation (\ref{Sec4import}) implies that the right hand side of (\ref{SCF}) can be written as
\begin{align*}
&
\sum_{(\alpha_1,\ldots, \alpha_q)\in C} \e V_{\alpha_1}' \cdots V_{\alpha_q}'  \prod_{\ell\leq q} \prod_{i\in C_\ell^1} \xi_i^{\alpha_\ell}
=
\sum_{(\rho(\alpha_1),\ldots, \rho(\alpha_q))\in C} \e V_{\rho(\alpha_1)}' \cdots V_{\rho(\alpha_q)}'  \prod_{\ell\leq q} \prod_{i\in C_\ell^1} \xi_i^{\rho(\alpha_\ell)}
\\
&=
\sum_{(\alpha_1,\ldots, \alpha_q)\in C} \e V_{\alpha_1} \cdots V_{\alpha_q}  \prod_{\ell\leq q} \prod_{i\in C_\ell^1} \txi_i^{\alpha_\ell}
=
\sum_{(\alpha_1,\ldots, \alpha_q)\in C} \e V_{\alpha_1} \cdots V_{\alpha_q} \,
\e \prod_{\ell\leq q} \prod_{i\in C_\ell^1} \txi_i^{\alpha_\ell}.
\end{align*}
This proves that
$$
\sum_{(\alpha_1,\ldots, \alpha_q)\in C} \e V_{\alpha_1}\cdots V_{\alpha_q} \,
\e \prod_{\ell\leq q} \prod_{i\in C_\ell^1} s_i^{\alpha_\ell}
=
\sum_{(\alpha_1,\ldots, \alpha_q)\in C} \e V_{\alpha_1} \cdots V_{\alpha_q} \,
\e \prod_{\ell\leq q} \prod_{i\in C_\ell^1} \txi_i^{\alpha_\ell}.
$$
Again, the second expectation in the sum on the right does not depend on $(\alpha_1,\ldots, \alpha_q)\in C$ and, since the choice of the constraints in the definition of $C$ was arbitrary, this proves that
\begin{equation}
\e \prod_{\ell\leq q} \prod_{i\in C_\ell^1} s_i^{\alpha_\ell}
=
\e \prod_{\ell\leq q} \prod_{i\in C_\ell^1} \txi_i^{\alpha_\ell}
\end{equation}
for any $\alpha_1,\ldots, \alpha_q \in \Natural$. Clearly, one can express any joint moment of the elements in these two arrays by choosing $q\geq 1$ large enough and choosing $\alpha_1,\ldots, \alpha_q$ and the sets $C_\ell^1$ properly, so the proof is complete.
\qed

\section{A consequence of the cavity equations}\label{Sec5label}

If we recall that  $\e_\alpha$ denotes the conditional expectation given the $\sigma$-algebra $\FF$ (i.e. in the random variables $\omega_{\alpha}$ and $\omega_{\alpha}^I$) then, using replicas and (\ref{Rab}), we can write for $\alpha, \beta, \gamma, \delta\in\Natural$ all different,
\begin{align}
&
\e \bigl(\e_\alpha s_1^{\alpha} s_2^{\alpha} - \e_\alpha s_1^{\alpha} \e_\alpha s_2^{\alpha} \bigr)^2
=
\e  s_1^{\alpha} s_2^{\alpha}  s_1^{\beta} s_2^{\beta} 
-2 \e  s_1^{\alpha} s_2^{\alpha}  s_1^{\beta} s_2^{\gamma}
+\e  s_1^{\alpha} s_1^{\beta}  s_2^{\gamma} s_2^{\delta}
\nonumber
\\
&= 
R_{\alpha,\beta}^2 - 2R_{\alpha,\beta}R_{\alpha,\gamma} + R_{\alpha,\beta}R_{\gamma,\delta} = (q_*)^2 - 2(q_*)^2 +(q_*)^2 = 0. 
\label{Sec5eqfirst}
\end{align}
By Theorem \ref{Sec4Th}, this implies that
$$
0 = \e  \txi_1^{\alpha} \txi_2^{\alpha}  \txi_1^{\beta} \txi_2^{\beta} 
-2 \e  \txi_1^{\alpha} \txi_2^{\alpha}  \txi_1^{\beta} \txi_2^{\gamma}
+\e  \txi_1^{\alpha} \txi_1^{\beta}  \txi_2^{\gamma} \txi_2^{\delta}
=
\e \bigl(\e_\alpha  \txi_1^{\alpha}  \txi_2^{\alpha} - \e_\alpha  \txi_1^{\alpha} \e_\alpha  \txi_2^{\alpha} \bigr)^2
$$
and, therefore,
$
\e_\alpha  \txi_1^{\alpha}  \txi_2^{\alpha} = \e_\alpha  \txi_1^{\alpha} \e_\alpha  \txi_2^{\alpha}
$
almost surely. If we recall that, conditionally on $\FF$, $(\txi_i^{\alpha})_{i\leq n}$ have the distribution of $(\xi_i^{\alpha})_{i\leq n}$ under the change of density (\ref{Ralpha}), we can rewrite this as
\begin{equation}
\e_\alpha  \xi_1^{\alpha}  \xi_2^{\alpha} R^\alpha = \e_\alpha  \xi_1^{\alpha} R^\alpha 
\, \e_\alpha  \xi_2^{\alpha} R^\alpha
\label{Sec5eq1}
\end{equation}
almost surely. Since this equation involves only two spin coordinates $i=1,2$, we can take $n=2$ in the definition of $R^\alpha$ as well, so that (recall (\ref{Aialpha}))
$$
\exp \zeta A^\alpha = \exp \zeta A_1^\alpha \, \exp \zeta A_2^\alpha.
$$
Let us denote by $\e_{\alpha,i}$ the expectation in the random variables $\omega_\alpha^{j,i,k}$ for $j\leq p-1$ and $k\geq 1$ that appear in the definition of $A_i^\alpha(\eps)$ in (\ref{Aieps}). Let us define 
\begin{align}
B_i^\alpha &=\, \frac{1}{\zeta}\log \e_{\alpha,i} \exp\zeta A_i^\alpha,
\\
Q^\alpha &=\, \frac{\exp\zeta(B_1^\alpha + B_2^\alpha)}{\e_\alpha \exp\zeta(B_1^\alpha + B_2^\alpha)},
\\
\eta_i^\alpha &=\,  \e_{\alpha,i} \xi_i^\alpha \exp \zeta(A_i^\alpha -B_i^\alpha).
\label{Sec5eqsecond}
\end{align}
Then it is easy to see that (\ref{Sec5eq1}) can be rewritten as
\begin{equation}
\e_\alpha  \eta_1^{\alpha}  \eta_2^{\alpha} Q^\alpha = \e_\alpha  \eta_1^{\alpha} Q^\alpha 
\, \e_\alpha  \eta_2^{\alpha} Q^\alpha
\label{Sec5eq2}
\end{equation}
almost surely. Since we already averaged the random variables $\omega_{\alpha}^I$, here the expectation $\e_\alpha$ is in $\omega_{\alpha}$ only. We will now use this to prove the main result of this section.
\begin{theorem}\label{Sec5Th}
The random variables $\eta_i^\alpha$ do not depend on $\omega_{\alpha}$.
\end{theorem}
Here and below, when we say that a function (or random variable) does not depend on a certain coordinate, this means that the function is equal to the average over that coordinate almost surely. Before we start the proof, let us make some simple preliminary observations. Both sides of (\ref{Sec5eq2}) depend on the random variables $g_{i,k}$, $\pi_i(\lambda p), \omega, \omega^{j,i,k}$ that generate the $\sigma$-algebra $\FF$. The Poisson random variables $\pi_1(\lambda p)$ and $\pi_2(\lambda p)$ can take any value $n\in \Natural$ at the same time with positive probability and, since (\ref{Sec5eq2}) holds almost surely, we can fix $\pi_1(\lambda p) = \pi_2(\lambda p) =n$ in (\ref{Sec5eq2}) for any $n\in \Natural$. Next, by the definition of the function $\theta$, both sides of (\ref{Sec5eq2}) are continuous functions of the variables $g_{i,k}$ for $k\leq n$. This implies that we can set $g_{k,1}=g_{k,2}=g_k$, and the equality (\ref{Sec5eq2}) will hold for all values of $g_k$, almost surely over $\omega$ and $\omega^{j,i,k}$ for $j\leq p-1$ and $k\leq n$. Finally, let us fix any $\omega$ such that (\ref{Sec5eq2}) holds for almost all $\omega^{j,i,k}$ for $j\leq p-1$ and $k\leq n$. Thus, from now on $\pi_1(\lambda p) = \pi_2(\lambda p) =n$, $g_{k,1}=g_{k,2}=g_k$ for $k\leq n$ and $\omega$ are all fixed. For simplicity of notation, let us temporarily denote $u_i = (\omega^{j,i,k})_{j\leq p-1, k\leq n}$ for $i=1,2$. Then 
$$
\eta_i^\alpha = \varphi(u_i,\omega_{\alpha}) \ \mbox{ and}\
Q^\alpha = \psi(u_1,u_2,\omega_{\alpha})
$$
for some functions $\varphi$ and $\psi$ and (\ref{Sec5eq2}) can be written as
\begin{equation}
\e_\alpha  \varphi(u_1,\omega_{\alpha}) \varphi(u_2,\omega_{\alpha})  \psi(u_1,u_2,\omega_{\alpha})
= 
\prod_{i=1,2}\e_\alpha  \varphi(u_i,\omega_{\alpha})  \psi(u_1,u_2,\omega_{\alpha})
\label{Sec5eq3}
\end{equation}
for almost all $u_1,u_2$. To prove Theorem \ref{Sec5Th}, we need to show that $\varphi(u,\omega_{\alpha})$ does not depend on $\omega_{\alpha}$.

\medskip
\noindent 
\textbf{Proof of Theorem \ref{Sec5Th}.}
If for fixed $\pi_1(\lambda p) = \pi_2(\lambda p) =n$ and $g_{k,1}=g_{k,2}=g_k$ for $k\leq n$ we denote $C = \beta\sum_{k\leq n}|g_k|$ then, by the definition of the function $\theta$, we can bound $A_i^\alpha(\eps)$ in (\ref{Aieps}) from above and below by $-C\leq A_i^\alpha(\eps)\leq C$. This implies that
\begin{equation}
e^{-4C}\leq Q^\alpha = \psi(u_1,u_2,\omega_{\alpha}) \leq e^{4C}.
\label{Sec5eq5}
\end{equation}
Of course, $|\varphi|\leq 1$. Suppose that for some $\eps>0$, there exists a set $U\subseteq [0,1]^{(p-1)n}$ of positive Lebesgue measure such that the variance $\mbox{Var}_{\omega_{\alpha}}(\varphi(u,\omega_{\alpha}))\geq \eps$ for $u\in U.$ Given $\delta>0$, let $(S_\ell)_{\ell\geq 1}$ be a partition of $L_1([0,1],dx)$ such that $\mbox{diam}(S_\ell)\leq \delta$ for all $\ell.$ Let
$$
U_\ell = \bigl\{ u\in [0,1]^{(p-1)n} \ |\ \varphi(u,\,\cdot\,) \in S_\ell \bigr\}.
$$
For some $\ell$, the Lebesgue measure of $U\cap U_\ell$ will be positive, so for some $u_1,u_2\in U$,
\begin{equation}
\e_{\alpha}| \varphi(u_1,\omega_{\alpha}) - \varphi(u_2,\omega_{\alpha}) | \leq \delta.
\label{Sec5eq6}
\end{equation}
The equations (\ref{Sec5eq5}) and (\ref{Sec5eq6}) imply that
$$
\bigl| 
\e_\alpha  \varphi(u_1,\omega_{\alpha}) \psi(u_1,u_2,\omega_{\alpha})
- \e_\alpha  \varphi(u_2,\omega_{\alpha}) \psi(u_1,u_2,\omega_{\alpha})
\bigr|
\leq e^{4C}\delta
$$
and, similarly,
$$
\bigl| 
\e_\alpha  \varphi(u_1,\omega_{\alpha}) \varphi(u_2,\omega_{\alpha})  \psi(u_1,u_2,\omega_{\alpha})
- \e_\alpha  \varphi(u_1,\omega_{\alpha})^2 \psi(u_1,u_2,\omega_{\alpha})
\bigr|
\leq e^{4C}\delta.
$$
Since $|\varphi|\leq 1$ and $\e_\alpha \psi =1$, the first inequality implies that
$$
\Bigl| 
\prod_{i=1,2}\e_\alpha  \varphi(u_i,\omega_{\alpha})  \psi(u_1,u_2,\omega_{\alpha})
- \bigl(\e_\alpha  \varphi(u_1,\omega_{\alpha}) \psi(u_1,u_2,\omega_{\alpha}) \bigr)^2
\Bigr|
\leq
e^{4C}\delta, 
$$
which, together with the second inequality and (\ref{Sec5eq3}), implies
$$
\e_\alpha  \varphi(u_1,\omega_{\alpha})^2 \psi(u_1,u_2,\omega_{\alpha})
-
 \bigl(\e_\alpha  \varphi(u_1,\omega_{\alpha}) \psi(u_1,u_2,\omega_{\alpha}) \bigr)^2
 \leq 2 e^{4C}\delta.
$$
The left hand side is a variance with the density $\psi$ and can be written using replicas as
$$
\frac{1}{2} \iint\! \bigl(\varphi(u_1,x)- \varphi(u_1,y)\bigr)^2 \psi(u_1,u_2,x)\psi(u_1,u_2,y)\,dx dy.
$$
By (\ref{Sec5eq5}) and the fact that $u_1\in U$, we can bound this from below by
$$
\frac{1}{2}e^{-8C} \iint\! \bigl(\varphi(u_1,x)- \varphi(u_1,y)\bigr)^2 \,dx dy
=
e^{-8C}\mbox{Var}_{\omega_{\alpha}}(\varphi(u_1,\omega_{\alpha})) \geq e^{-8C}\eps.
$$
Comparing lower and upper bounds, $e^{-8C}\eps \leq e^{4C}\delta$, we arrive at contradiction, since $\delta>0$ was arbitrary. Therefore, $\mbox{Var}_{\omega_{\alpha}}(\varphi(u,\omega_{\alpha})) = 0$ for almost all $u$ and this finishes the proof.
\qed

\section{Proof of Theorem \ref{Th1} when $h=0$} \label{Sec6label}

We begin with one basic observation. For integer $m\geq 1$, let us define 
\begin{equation}
f^{(m)}(w,u,v) = \int_0^1 \! f(w,u,v,x)^m \, dx.
\label{fbareq}
\end{equation}
The equation (\ref{Sec3eq1}) implies that
\begin{equation}
\int_0^1 \! f^{(1)}(w, u_1, v) f^{(1)}(w,u_2, v) \, dv
= q_*
\label{Sec4eq2}
\end{equation}
for almost all $w,u_1,u_2\in [0,1]$, which in turn implies the following.
\begin{lemma} \label{Sec4Lem1}
The function $f^{(1)}(w,u,v)$ does not depend on the second coordinate $u$.
\end{lemma}
\textbf{Proof.} By (\ref{Sec4eq2}), for almost all $w\in[0,1]$, the  measure 
$$
du\circ \bigl(u\to f^{(1)}(w,u,\cdot) \bigr)^{-1}
$$ 
on $L^2([0,1],dv)$ (or $H$ in (\ref{spaceH})) is such that, for any two points $\sigma^1,\sigma^2$ in its support, we have $\sigma^1\cdot\sigma^2 = q_*.$ Clearly, this can happen only if the measure is concentrated on one point and this finishes the proof.
\qed

\medskip
\medskip
\noindent
For simplicity of notation, we will sometimes omit the coordinate $u$ but still use the same notation for the function. For example, we will write $f^{(1)}(w,v)$ and notice that, by
Lemma \ref{Sec4Lem1} and (\ref{Sec4eq2}), 
\begin{equation}
\int_0^1 \! f^{(1)}(w, v)^2 \, dv
= q_*
\label{Sec6qstar}
\end{equation}
for almost all $w\in [0,1]$. This implies another observation, which requires no proof.
\begin{lemma}\label{Sec6Lem1}
If $q_*\not =0$ then, for almost all $(w,u)\in [0,1]^2$, the functions $f^{(1)}(w,u, \,\cdot\, )=f^{(1)}(w, \,\cdot\,)$ and $f^{(m)}(w,u,\,\cdot\,)$ for even $m\geq 2$ are not identically zero.
\end{lemma}
In Theorem \ref{Sec5Th} we proved that $\eta_i^\alpha$ does not depend on $\omega_{\alpha}$ and, tracing all the definitions back to (\ref{Aieps}), $\eta_i^\alpha$ can be written as 
\begin{equation}
\eta_i^\alpha 
=
 \frac{\e_{\alpha, i} \Av \eps \exp A_i^\alpha(\eps) (\Av \exp A_i^\alpha(\eps))^{\zeta-1}}{\e_{\alpha,i} (\Av \exp A_i^\alpha(\eps))^\zeta}.
\label{eta}
\end{equation}
Since $\eps\in \{-1,+1\}$, we can write 
$$
\exp h\eps = \ch(h)\bigl(1+\myth(h)\eps\bigr).
$$
Therefore, by (\ref{DefthetaG}) and (\ref{Aieps}),
$$
\exp A_i^\alpha(\eps)
=
\ch(h)
\prod_{k\leq \pi_i(\lambda p)} \ch(\beta g_{i,k})\bigl(1+ \myth(\beta g_{i,k})\, s_{1,i,k}^\alpha\cdots s_{p-1,i,k}^\alpha \eps \bigr) \bigl(1+\myth(h)\eps\bigr).
$$
Obviously, $\ch(h)$ and all the factors $\ch(\beta g_{i,k})$ will cancel out in (\ref{eta}) so we can omit them. Let us now fix $\pi_i(\lambda p) = n$. Since $\alpha$ and $i$ are now fixed, for simplicity of notation, let us denote 
\begin{equation}
t_k = \myth(\beta g_{i,k})
\ \mbox{ and }\
s_{(k)}=s_{1,i,k}^\alpha\cdots s_{p-1,i,k}^\alpha
\label{Sec6eq1}
\end{equation}
and for the remainder of the paper redefine 
\begin{equation}
\exp A_i^\alpha(\eps)
=
\prod_{k\leq n} \bigl(1+ t_k\, s_{(k)} \eps \bigr)\bigl(1+\myth(h)\eps\bigr).
\label{AiSec6}
\end{equation}
As in the discussion in the paragraph below Theorem \ref{Sec5Th}, since $\pi_i(\lambda p)$ is discrete and $\eta_i^\alpha$ is a continuous function of $t_k = \myth(\beta g_{i,k})$, we can say that  $\eta_i^\alpha$ does not depend on $\omega_\alpha$ in the sense that it is equal to its expectation $\e_{\alpha}$ in $\omega_\alpha$ for all $n\geq 1$, for almost all $\omega_\alpha$, for almost all $u_i = (\omega^{j,i,k})_{j\leq p-1, k\leq n}$ and for all $t_k\in (-1,1)$ for $k\leq n$. In particular, we will use that, for all $n\geq 1$, all partial derivatives of $\eta_i^\alpha$ in $(t_k)_{k\leq n}$ do not depend on $\omega_\alpha$. 

In the remainder of this section we will consider the more difficult case when $h=0$. If we take $n=3$ then from (\ref{AiSec6}) we obtain
\begin{align}
&
X:= \Av \exp A_i^\alpha(\eps) = 1+ t_1 t_2 s_{(1)} s_{(2)}+ t_1 t_3 s_{(1)} s_{(3)}+ t_2 t_3 s_{(2)} s_{(3)},
\nonumber
\\
&
Y:= \Av \eps \exp A_i^\alpha(\eps) = t_1 s_{(1)} +t_2 s_{(2)} +t_3 s_{(3)} + t_1 t_2 t_3 s_{(1)}s_{(2)} s_{(3)}.
\label{Sec6XY}
\end{align}
We will fix $t_3$ to be any non-zero value, for example, $t_3=1/2$ and, for any $m\geq 1$, consider
\begin{align}
\frac{\partial^m \eta_i^\alpha}{\partial t_1 \partial t_2^{m-1}}\,\Bigr|_{t_1=t_2=0}
&=\,
\frac{\partial^m}{\partial t_1 \partial t_2^{m-1}}
 \frac{\e_{\alpha, i} Y X^{\zeta-1}}{\e_{\alpha,i} X^\zeta}
\,\Bigr|_{t_1=t_2=0},
\label{deriv}
\\
\frac{\partial^m \eta_i^\alpha}{\partial t_2^{m}}\,\Bigr|_{t_1=t_2=0}
&=\,
\frac{\partial^m}{\partial t_2^{m}}
 \frac{\e_{\alpha, i} Y X^{\zeta-1}}{\e_{\alpha,i} X^\zeta}
\,\Bigr|_{t_1=t_2=0}.
\label{deriv2}
\end{align}
We will prove the following.
\begin{lemma} \label{Sec6Lem2}
The derivative (\ref{deriv}) is given by a linear combination of $\e_{\alpha,i} s_{(1)}^{{}}s_{(2)}^{m-1} s_{(3)}^{m+1}$ with some non-zero coefficient and various products of factors of the type 
\begin{equation}
\e_{\alpha,i} s_{(1)}^{m_1}s_{(2)}^{m_2} s_{(3)}^{m_3}
\label{Sec6factor}
\end{equation}
with integer powers $m_1,m_2,m_3\leq m$. The derivative (\ref{deriv2}) is given by a linear combination of $\e_{\alpha,i} s_{(2)}^{m} s_{(3)}^{m+1}$ with some non-zero coefficient and various products of factors of the type (\ref{Sec6factor}) with integer powers $m_1,m_2,m_3\leq m$.
\end{lemma}
\textbf{Proof.}
We will only prove the first claim concerning the derivative (\ref{deriv}), since the proof of the second claim is similar. First of all, notice that when we apply a derivative to the denominator, this results in some power of the denominator and we get another factor equal to the derivative of the $\e_{\alpha,i} X^\zeta$. Since $X|_{t_1=t_2=0} =1$, in the end all denominators will just be equal to one. When we take a derivative of some power of $X$, we end up with another power of $X$ times one of the factors
\begin{equation}
\frac{\partial X}{\partial t_1} = t_2 s_{(1)} s_{(2)}+ t_3 s_{(1)} s_{(3)},\,\,
\frac{\partial X}{\partial t_2} = t_1 s_{(1)} s_{(2)}+t_3 s_{(2)} s_{(3)}.
\label{Sec6factors}
\end{equation}
Further non-trivial derivatives of such factors can only give us
$$
\frac{\partial^2 X}{\partial t_1 \partial t_2} = s_{(1)} s_{(2)}.
$$
On the other hand, if these factors are not differentiated, in the end they become
$$
\frac{\partial X}{\partial t_1}\,\Bigr|_{t_1=t_2=0} =  t_3 s_{(1)} s_{(3)},\,\,
\frac{\partial X}{\partial t_2}\,\Bigr|_{t_1=t_2=0} = t_3 s_{(2)} s_{(3)}.
$$
Non-trivial derivatives of $Y$ will include $Y|_{t_1=t_2=0} =t_3 s_{(3)}$ and
\begin{equation}
\frac{\partial Y}{\partial t_1}\,\Bigr|_{t_1=t_2=0} = s_{(1)},\,\,
\frac{\partial Y}{\partial t_2}\,\Bigr|_{t_1=t_2=0} = s_{(2)}
\ \mbox{ and }\
\frac{\partial^2 Y}{\partial t_1 \partial t_2}\,\Bigr|_{t_1=t_2=0} = t_3 s_{(1)}s_{(2)} s_{(3)}.
\label{Sec6dY}
\end{equation}
All together, this makes it clear that the derivative (\ref{deriv}) will be given by a sum of various products of factors of the type $\e_{\alpha,i} s_{(1)}^{m_1}s_{(2)}^{m_2} s_{(3)}^{m_3}$. 

The main observation we will need is the answer to the following question: what is the largest power among $m_1, m_2$ and $m_3$ that we can possibly achieve? Let us first present one candidate for the answer. In (\ref{deriv}), let us not touch the denominator, let us not differentiate the factor $Y$, and apply all derivatives to the factor $X^{\zeta-1}$. In the end, the factor $Y$ will give $Y|_{t_1=t_2=0} =t_3 s_{(3)}$. Also, every time we differentiate the power of $X$ and get one of the factors in (\ref{Sec6factors}), let us not differentiate these factors and continue differentiating only the powers of $X$. We will end up with the term (up to a non-zero constant)
\begin{equation}
\frac{1}{\e_{\alpha,i} X^\zeta} \e_{\alpha, i} \Bigl(Y X^{\zeta-1-m}
\frac{\partial X}{\partial t_1}\Bigl(\frac{\partial X}{\partial t_2}\Bigr)^{m-1} \Bigr)
\,\Bigr|_{t_1=t_2=0}
=
t_3^{m+1} \e_{\alpha,i} s_{(1)}^{{}}s_{(2)}^{m-1} s_{(3)}^{m+1}. 
\label{Sec6max}
\end{equation}
Notice that we get one power of $s_{(3)}$ from the factor $Y$ and each time we differentiate a power of $X$, we gain one power of $s_{(3)}$. It remains to understand why there is no other way to obtain the power $m+1$ for one of the factors. The key point is that, for $a=\zeta$ or $a=\zeta-1$ and any $k\geq 1$, the derivatives
$$
\frac{\partial^k X^a}{\partial t_1 \partial t_2^{k-1}}\,\Bigr|_{t_1=t_2=0}
\ \mbox{ and }\
\frac{\partial^k X^a}{\partial t_2^{k}}\,\Bigr|_{t_1=t_2=0}
$$  
can not produce a power higher than $s_{(1)}^k$ or $s_{(2)}^k$ or $s_{(3)}^k$. This can be proved formally by considering the Taylor series for $(1+x)^a$ around $x=0$ with
$$
x=t_1 t_2 s_{(1)} s_{(2)}+ t_1 t_3 s_{(1)} s_{(3)}+ t_2 t_3 s_{(2)} s_{(3)}.
$$ 
There are two ways to get a term with $t_1 t_2^{k-1}$,
\begin{equation}
t_1 t_2 s_{(1)} s_{(2)} \bigl(t_2 t_3 s_{(2)} s_{(3)} \bigr)^{k-1}
\ \mbox{ or }\
t_1 t_3 s_{(1)} s_{(3)}\bigl(t_2 t_3 s_{(2)} s_{(3)} \bigr)^{k-1},
\label{Sec6term1}
\end{equation}
and there is only one way to get a term with $t_2^{k}$ and without $t_1$, $(t_2 t_3 s_{(2)} s_{(3)})^{k}$. In both cases, the highest power is $k$. Let us now consider various cases. If $k=m$, we are computing the derivative $\partial^m/\partial t_1\partial t_2^{m-1}$, so the terms (\ref{Sec6term1}) will give us (up to constants)
$$
s_{(1)} s_{(2)}^m s_{(3)}^{m-1}
\ \mbox{ or }\
s_{(1)} s_{(2)}^{m-1} s_{(3)}^{m}.
$$
If we apply this derivative to $X^{\zeta-1}$, with another factor $s_{(3)}$ coming from $Y|_{t_1=t_2=0}$, we will get 
$$
s_{(1)} s_{(2)}^m s_{(3)}^{m}
\ \mbox{ or }\
s_{(1)} s_{(2)}^{m-1} s_{(3)}^{m+1}.
$$
Of course, the second term is the one we got in (\ref{Sec6max}). Now, let us consider other possibilities that do not involve differentiating the denominator. If we waste one or two derivatives on $Y$ as in (\ref{Sec6dY}), all the factors  still have power one, and now we are left with at most $k=m-1$ derivatives to apply to $X^{\zeta-1}$. We know that the highest power we can achieve is $k=m-1$, and it is not enough to reach $m+1$. Finally, if we apply some derivative to the denominator, the best we can do is to apply all derivatives to the factor $X^\zeta$. In this case, again, we can only attain the highest power equal to $m$ and this proves the first claim. For the second claim, the proof is almost the same using
\begin{equation}
\frac{1}{\e_{\alpha,i} X^\zeta} \e_{\alpha, i} \Bigl(Y X^{\zeta-1-m}
\Bigl(\frac{\partial X}{\partial t_2}\Bigr)^{m} \Bigr)
\,\Bigr|_{t_1=t_2=0}
=
t_3^{m+1} \e_{\alpha,i} s_{(2)}^{m} s_{(3)}^{m+1}
\label{Sec6max2}
\end{equation}
instead of (\ref{Sec6max}). This finishes the proof.
\qed

\medskip
\noindent
Let us recall the function $f^{(m)}(w,u,v)$ defined in (\ref{fbareq}) and consider a function
\begin{equation}
g^{(m)}\bigl(w, u, (v_j)_{j\leq p-1}\bigr) := 
\prod_{j\leq p-1}f^{(m)}(w,u,v_{j}).
\label{Sec7fprod}
\end{equation} 
We will now show that Lemma \ref{Sec6Lem2} implies the following.
\begin{lemma} \label{Sec6Lem3}
For all $m\geq 1$, the function $g^{(m)}$ defined in (\ref{Sec7fprod}) does not depend on $u$ and, therefore, we can write it as $g^{(m)}\bigl(w, (v_j)_{j\leq p-1}\bigr)$.
\end{lemma}
\textbf{Proof.} For better agreement with the notation in Lemma \ref{Sec6Lem2}, we will be proving that $g^{(m+1)}$ does not depend on $u$ for $m\geq 0$. By Lemma \ref{Sec4Lem1}, we know this for $m=0$. For $m=1$, the derivative (\ref{deriv}) will be a linear combination of $\e_{\alpha,i} s_{(1)}^{{}}s_{(3)}^{2}$ and other terms consisting of products of factors in (\ref{Sec6factor}) with powers $m_k$ equal to $0$ or $1$. 

Since the expectation $\e_{\alpha,i}$ is in the random variables $\omega_\alpha^{j,i,k}$ and these random variables are independent in different factors $s_{(k)}$, we get
$$
\e_{\alpha,i} s_{(1)}^{m_1}s_{(2)}^{m_2} s_{(3)}^{m_3}
=
\e_{\alpha,i} s_{(1)}^{m_1} \e_{\alpha,i} s_{(2)}^{m_2} \e_{\alpha,i}  s_{(3)}^{m_3}.
$$
Furthermore, for the same reason
$$
\e_{\alpha,i}s_{(k)}=\e_{\alpha,i}s_{1,i,k}^\alpha\cdots \e_{\alpha,i}s_{p-1,i,k}^\alpha.
$$
Finally, by (\ref{fbareq}) and Lemma \ref{Sec4Lem1},
$$
\e_{\alpha,i} s_{j,i,k}^\alpha
=
\e_{\alpha,i} f(\omega, \omega_\alpha, \omega^{j,i,k},\omega_\alpha^{j,i,k})
=
f^{(1)}(\omega, \omega_\alpha, \omega^{j,i,k})
=
f^{(1)}(\omega, \omega^{j,i,k})
$$
almost surely, so it does not depend on $\omega_\alpha$. This proves that all the factors (\ref{Sec6factor}) with powers $m_k$ equal to $0$ or $1$ do not depend on $\omega_\alpha$ and, since the derivative (\ref{deriv}) also does not depend on $\omega_\alpha$, 
\begin{equation}
\e_{\alpha,i} s_{(1)}^{{}}s_{(3)}^{2}
=
\e_{\alpha,i} s_{(1)}^{{}} \e_{\alpha,i} s_{(3)}^{2}
=
\e_{\alpha,i} s_{(3)}^{2}
\prod_{j\leq p-1} f^{(1)}(\omega, \omega^{j,i,1})
\label{Sec6prod}
\end{equation}
does not depend on $\omega_\alpha$. When $q_*=0$, by (\ref{Sec6qstar}), the right hand side in (\ref{Sec6prod}) is equal to zero almost surely and we get no information. When $q_*\not = 0$, by Lemma \ref{Sec6Lem1}, for almost all $\omega$, we can find $(\omega^{j,i,1})_{j\leq p-1}$ such that the last product on the right hand side of (\ref{Sec6prod}) is not zero and, since it does not depend on $\omega_\alpha$, we proved that 
\begin{equation}
\e_{\alpha,i} s_{(3)}^{2} = \prod_{j\leq p-1} f^{(2)}(\omega, \omega_\alpha, \omega^{j,i,3})
\label{Sec6there}
\end{equation}
does not depend on $\omega_\alpha$. In other words, the function (\ref{Sec7fprod}) for $m=2$ does not depend on $u$.

Suppose that we proved that the function $g^{(\ell)}$ does not depend on $u$ for $\ell\leq m$. To make the induction step, we will use the first statement of Lemma \ref{Sec6Lem2} for odd $m$, and the second statement for even $m$. In both cases, by the induction assumption, each factor in (\ref{Sec6factor}),
$$
\e_{\alpha,i} s_{(1)}^{m_1}s_{(2)}^{m_2} s_{(3)}^{m_3}
=
\prod_{k\leq 3} \e_{\alpha,i} s_{(k)}^{m_k}
=
\prod_{k\leq 3} \prod_{j\leq p-1}
\e_{\alpha,i}(s_{j,i,k}^\alpha)^{m_k}
=
\prod_{k\leq 3} \prod_{j\leq p-1}
f^{(m_k)}(\omega, \omega_\alpha, \omega^{j,i,k}),
$$
does not depend on $\omega_\alpha$, because all $m_k\leq m$. Since the derivatives (\ref{deriv}) and (\ref{deriv2}) do not depend on $\omega_\alpha$, Lemma \ref{Sec6Lem2} implies that
$$
\e_{\alpha,i} s_{(1)}^{{}}s_{(2)}^{m-1} s_{(3)}^{m+1}
\ \mbox{ and }\
\e_{\alpha,i} s_{(2)}^{m} s_{(3)}^{m+1}
$$
do not depend on $\omega_\alpha$. 

When $m$ is odd, we use that $\e_{\alpha,i} s_{(1)}^{{}}s_{(2)}^{m-1} s_{(3)}^{m+1}$, which by the induction assumption is equal to
\begin{align*}
&
\prod_{j\leq p-1} f^{(1)}(\omega, \omega_\alpha, \omega^{j,i,1})
\prod_{j\leq p-1} f^{(m-1)}(\omega, \omega_\alpha, \omega^{j,i,2})
\prod_{j\leq p-1} f^{(m+1)}(\omega, \omega_\alpha, \omega^{j,i,3})
\\
& =
g^{(1)}\bigl(\omega, (\omega^{j,i,1})_{j\leq p-1}\bigr)
g^{(m-1)}\bigl(\omega, (\omega^{j,i,2})_{j\leq p-1} \bigr)
\prod_{j\leq p-1} f^{(m+1)}(\omega, \omega_\alpha, \omega^{j,i,3}),
\end{align*}
does not depend on $\omega_\alpha$. Because $m-1$ is even in this case, by Lemma \ref{Sec6Lem1}, for almost all $\omega$ we can find $(\omega^{j,i,1})_{j\leq p-1}$ and $(\omega^{j,i,2})_{j\leq p-1}$ such that the first two factors are not zero, and this implies that
\begin{equation}
\e_{\alpha,i} s_{(3)}^{m+1} = 
\prod_{j\leq p-1} f^{(m+1)}(\omega, \omega_\alpha, \omega^{j,i,3})
\label{Sec6there2}
\end{equation}
does not depend on $\omega_\alpha$. This completes the induction step when $m$ is odd.

When $m$ is even, we use that $\e_{\alpha,i} s_{(2)}^{m} s_{(3)}^{m+1}$, which by the induction assumption is equal to
\begin{align*}
&
\prod_{j\leq p-1} f^{(m)}(\omega, \omega_\alpha, \omega^{j,i,2})
\prod_{j\leq p-1} f^{(m+1)}(\omega, \omega_\alpha, \omega^{j,i,3})
\\
& =
g^{(m)}\bigl(\omega, (\omega^{j,i,2})_{j\leq p-1} \bigr)
\prod_{j\leq p-1} f^{(m+1)}(\omega, \omega_\alpha, \omega^{j,i,3}),
\end{align*}
does not depend on $\omega_\alpha$. Because $m$ is even, by Lemma \ref{Sec6Lem1}, for almost all $\omega$ we can find $(\omega^{j,i,2})_{j\leq p-1}$ such that the first factor is not zero, and this again implies that
\begin{equation}
\e_{\alpha,i} s_{(3)}^{m+1} = 
\prod_{j\leq p-1} f^{(m+1)}(\omega, \omega_\alpha, \omega^{j,i,3})
\end{equation}
does not depend on $\omega_\alpha$. This finishes the proof.
\qed

\medskip
\noindent
We are now ready to prove the first claim of Theorem \ref{Th1}.

\medskip
\noindent
\textbf{Proof of Theorem \ref{Th1}} (\emph{The case $h=0$}). Let us consider functions $T_j: [0,1]\to [0,1]$ for $j\leq p-1$ such that $(T_j(\omega))_{j\leq p-1}$ are i.i.d. uniform on $[0,1]$ when $\omega$ is uniform on $[0,1]$. Consider a function $g:[0,1]^4 \to [-1,1]$ given by
\begin{equation}
g(w,u,v,x) = \prod_{j\leq p-1} f\bigl(w,u,T_j(v),T_j(x) \bigr).
\end{equation}
Then, for any $m\geq 1$,
$$
g^{(m)}(w,u,v):=
\int\! g(w,u,v,x)^m \, dx
=
\prod_{j\leq p-1}
\int\! f(w,u,T_j(v),x)^m \, dx
=
\prod_{j\leq p-1} f^{(m)}(w,u,T_j(v)).
$$
Since we showed in Lemma \ref{Sec6Lem3} that the right hand side does not depend on $u$, if we consider the conditional distribution of $g(w,u,v,x)$ given $(w,u,v)$,
\begin{equation}
\p(w,u,v; [-\infty,y]) = \bigl|\bigl\{x \ \bigr|\  g(w,u,v,x)\leq y\bigr\}\bigr|
\label{Sec6distP}
\end{equation}
then this distribution function does not depend on $u$ and we can write it as $\p(w,v; [-\infty,y])$. If we consider its quantile transform
$$
h(w,v,x) = \inf\bigl\{y \ \bigr|\ x\leq \p(w,v; [-\infty,y])\bigr\}
$$
then, for all $m\geq 1$, we have
$$
\int_0^1 \! g(w,u,v,x)^m \, dx
=
\int_0^1 \! h(w,v,x)^m \, dx
$$
almost surely. By comparing the joint moments, this implies that
\begin{equation}
\bigl( g(\omega, \omega_\alpha, \omega^I, \omega_\alpha^I) \bigr)_{\alpha\in \Natural, I\in \II }
\, \stackrel{d}{=}\,
\bigl( h(\omega, \omega^I, \omega_\alpha^I) \bigr)_{\alpha\in \Natural, i\in\II}
\label{sigmaf2end}
\end{equation}
for any countable set of multi-indices $\II$. On the other hand, if we take $I=(i,k)$, then
$$
g\bigl(\omega, \omega_\alpha, \omega^{i,k}, \omega_\alpha^{i,k} \bigr)
=
\prod_{j\leq p-1} f\bigl(\omega,\omega_\alpha,T_j(\omega^{i,k}),T_j(\omega_\alpha^{i,k})\bigr)
$$
can be viewed, by the definition of the functions $(T_j)_{j\leq p-1}$, as another representation for
$$
z_{i,k}^\alpha : =
\prod_{j\leq p-1} s_{j,i,k}^\alpha
=
\prod_{j\leq p-1} f(\omega,\omega_\alpha,\omega^{j,i,k},\omega_\alpha^{j,i,k}).
$$
More specifically, the equation (\ref{sigmaf2end}) implies that
\begin{equation}
\bigl( z_{i,k}^\alpha \bigr)_{\alpha,i,k\in \Natural}
\, \stackrel{d}{=}\,
\bigl( h(\omega, \omega^{i,k}, \omega_\alpha^{i,k}) \bigr)_{\alpha,i,k\in \Natural}.
\label{sigmaf2end2}
\end{equation}
Since the cavity fields $A_i^{\alpha}(\eps)$ in (\ref{Aieps}) can be written as
$$
A_i^\alpha(\eps) 
= \sum_{k\leq \pi_i(\lambda p)} \beta g_{i,k} \prod_{j\leq p-1} s_{j,i,k}^\alpha \,\eps +h\eps
= \sum_{k\leq \pi_i(\lambda p)} \beta g_{i,k} z_{i,k}^\alpha \,\eps +h\eps,
$$
we can now redefine them in the cavity equations by directly setting
$$
z_{i,k}^\alpha 
=
h(\omega, \omega^{i,k}, \omega_\alpha^{i,k})
$$
instead of defining the factors $s_{j,i,k}^\alpha$ in (\ref{sialpha}) separately. Since $A_i^\alpha(\eps)$ now does not depend on $\omega_\alpha$, the change of density in (\ref{Ralpha}) can be rewritten as $R^\alpha = \prod_{i\leq n} R^\alpha_i$, where
\begin{equation}
R^\alpha_i : = \frac{\exp \zeta A_i^\alpha}{ \e_{\alpha,i} \exp \zeta A_i^\alpha}
\label{Ralphai}
\end{equation}
and where $\e_{\alpha,i}$ denotes the expectation in the random variables $(\omega_{\alpha}^{i,k})_{k\geq 1}$. Notice that both random variables $\xi_i^{\alpha}$ in (\ref{Sec3xiialpha}) and $R_{i}^\alpha$ are now functions of 
$$
\omega,\,\, U_i : = \bigl( \pi_i(\lambda p), (g_{i,k})_{k\geq 1}, (\omega^{i,k})_{k\geq 1}\bigr)
\ \mbox{ and }\
U_i^\alpha : = (\omega_\alpha^{i,k})_{k\geq 1}.
$$
The random variables $U_i$ are i.i.d. for $i\geq 1$, and the random variables $U_i^\alpha$ are i.i.d. for $\alpha, i\geq 1$. Therefore, since the change of density $R^\alpha$ decouples as in (\ref{Ralphai}), conditionally on the $\sigma$-algebra $\FF$ generated by the random variables $\omega$ and $(U_i)_{i\geq 1}$, the distribution of each $\txi_i^{\alpha}$ in Theorem \ref{Sec4Th} is now simply the distribution of $\xi_i^{\alpha}$ under the change of density $R_{i}^\alpha$. Moreover, conditionally on $\FF$, the random variables $\txi_i^{\alpha}$ are independent and can be generated in distribution as $\txi_i^{\alpha} = F(\omega,U_i, \omega_\alpha^i)$ for some function $F$. We can also generate $U_i$ as a function of $\omega^i$ uniform on $[0,1]$ and, therefore, in distribution, $\txi_i^{\alpha} = f(\omega,\omega^i, \omega_\alpha^i)$ for some function $f$. By Theorem \ref{Sec4Th} and (\ref{sigmaf2}), this precisely gives the representation (\ref{sigmaf}), so the proof if finished.
\qed

\section{Proof of Theorem \ref{Th1comp}}\label{Sec7label}

We will now prove Theorem \ref{Th1comp}. First, suppose that $q_*=0$. Then, similarly to (\ref{Sec5eqfirst}), we can write
\begin{align}
&
\e \bigl(\e_\alpha s_1^{\alpha} s_2^{\alpha} \bigr)^2
=
\e  s_1^{\alpha} s_2^{\alpha}  s_1^{\beta} s_2^{\beta} 
= 
R_{\alpha,\beta}^2 = (q_*)^2 = 0.
\label{Sec7eqfirst}
\end{align}
By Theorem \ref{Sec4Th}, this implies that
$$
0 = \e  \txi_1^{\alpha} \txi_2^{\alpha}  \txi_1^{\beta} \txi_2^{\beta} 
=
\e \bigl(\e_\alpha  \txi_1^{\alpha}  \txi_2^{\alpha}\bigr)^2
$$
and, therefore, $\e_\alpha  \txi_1^{\alpha}  \txi_2^{\alpha} = 0$ almost surely. In the notation (\ref{Sec5eqsecond}), this can be rewritten as
\begin{equation}
\e_\alpha  \eta_1^{\alpha}  \eta_2^{\alpha} Q^\alpha = 0
\label{Sec7eq2}
\end{equation}
almost surely. As in the proof of Theorem \ref{Sec5Th}, this implies that $\eta_i^\alpha=0$ almost surely or, equivalently,
\begin{equation}
\e_{\alpha, i} \Av \eps \exp A_i^\alpha(\eps) (\Av \exp A_i^\alpha(\eps))^{\zeta-1}=0
\label{eta2}
\end{equation}
almost surely. Recalling the notation (\ref{Sec6eq1}) and (\ref{AiSec6}), if we now take $n=2$ then we get
$$
\e_{\alpha, i} \bigl(t_1 s_{(1)}+t_2 s_{(2)} \bigr) \bigl(1+t_1 t_2 s_{(1)} s_{(2)} \bigr)^{\zeta-1} = 0
$$
almost surely. Using the Taylor series for $(1+x)^{\zeta-1}$, we see that the monomial $t_1^{m+1}  t_2^m$ appears with the factor 
\begin{equation}
\e_{\alpha, i} s_{(1)}^{m+1} s_{(2)}^m
=
\prod_{j\leq p-1} f^{(m+1)}(\omega, \omega_\alpha, \omega^{j,i,1})
\prod_{j\leq p-1} f^{(m)}(\omega, \omega_\alpha, \omega^{j,i,2}),
\label{Sec7eq3}
\end{equation}
which then must be equal to zero almost surely. Similarly to (\ref{Sec3eq1}), we can write
\begin{equation}
0\not = q^* = R_{\alpha,\alpha}= \e_i f(\omega, \omega_\alpha, \omega^i, \omega_\alpha^i)^2 
= \e_i f^{(2)}(\omega, \omega_\alpha, \omega^i). 
\label{Sec7Raa}
\end{equation}
This implies that for almost all $\omega, \omega_\alpha$ the function $f^{(2)}(\omega, \omega_\alpha, \,\cdot\,)$ is not identically zero which, of course, implies that $f^{(m)}(\omega, \omega_\alpha, \,\cdot\,)$ is not identically zero for all even $m\geq 2$. Therefore, for even $m$, we can find $(\omega^{j,i,2})_{j\leq p-1}$ such that the second product in (\ref{Sec7eq3}) is not zero, so
$$
\prod_{j\leq p-1} f^{(m+1)}(\omega, \omega_\alpha, \omega^{j,i,1}) = 0
$$
almost surely. Clearly, this implies that $f^{(m+1)}(w, u, v) = 0$ almost surely for all even $m$, which means that the conditional distribution of $f(w,u,v,x)$ given $(w,u,v)$,
\begin{equation}
\p(w,u,v; [-\infty,y]) = \bigl|\bigl\{x \ \bigr|\  f(w,u,v,x)\leq y\bigr\}\bigr|,
\label{Sec7distP}
\end{equation}
is symmetric.

Now, suppose that the distribution in (\ref{Sec7distP}) is symmetric for almost all $(w,u,v)$. Since, for any $n\geq 1$,
\begin{align*}
\Av \prod_{k\leq n}\bigl(1+ t_k (-s_{(k)}) \eps\bigr)
&=\,
\Av \prod_{k\leq n}\bigl(1+ t_k s_{(k)} \eps\bigr),
\\
\Av \eps \prod_{k\leq n}\bigl(1+ t_k (-s_{(k)}) \eps\bigr)
&=\,
- \Av \eps \prod_{k\leq n}\bigl(1+ t_k s_{(k)} \eps\bigr),
\end{align*}
the symmetry of the distribution $\p(w,u,v; [-\infty,y])$ implies (\ref{eta2}) and, thus, (\ref{Sec7eq2}). In turn, (\ref{Sec7eq2}) implies (\ref{Sec7eqfirst}), so $q_* = 0$.
\qed

\section{Proof of Theorem \ref{Th1}, the case $h\not= 0$} \label{Sec8label}  

For simplicity of notation, we will denote 
$
c=\myth(h)\in (-1,1)\setminus\{0\}.
$ 
If we take $n=1$ in (\ref{AiSec6}) then
\begin{equation}
\Av \exp A_i^\alpha(\eps) = 1+ c t_1 s_{(1)}
\ \mbox{ and }\
\Av \eps \exp A_i^\alpha(\eps) = c+ t_1 s_{(1)}.
\label{Sec6XYh}
\end{equation}
Now, using that both
$$
c^{-1}\eta_i^\alpha 
=
\frac{\e_{\alpha, i} (1+ c^{-1}t_1 s_{(1)}) (1+ c t_1 s_{(1)})^{\zeta-1}}{\e_{\alpha,i} (1+ c t_1 s_{(1)})^\zeta}
\ \mbox{ and }\
1 = \frac{\e_{\alpha, i} (1+ c t_1 s_{(1)}) (1+ c t_1 s_{(1)})^{\zeta-1}}{\e_{\alpha,i} (1+ c t_1 s_{(1)})^\zeta}
$$
do not depend on $\omega_\alpha$, we get that
\begin{equation}
\frac{1-c^{-1}\eta_i^\alpha}{c-c^{-1}}=
\frac{\e_{\alpha, i} t_1 s_{(1)} (1+ c t_1 s_{(1)})^{\zeta-1}}{\e_{\alpha,i} (1+ c t_1 s_{(1)})^\zeta}
\label{Sec8eq1}
\end{equation}
does not depend on $\omega_\alpha$. As in the proof of the case $h=0$ (only much easier) one can show that, for $m\geq 1$, the derivative $\partial^{m}/\partial t_1^{m}$ of the right hand side at $t_1=0$ will be a linear combination of $\e_{\alpha, i} s_{(1)}^{m}$ and various products of factors $\e_{\alpha, i} s_{(1)}^{m_1}$ for $m_1<m$. By induction on $m$, this implies that  
$$
\e_{\alpha,i} s_{(1)}^{m} = 
\prod_{j\leq p-1} f^{(m)}(\omega, \omega_\alpha, \omega^{j,i,1})
=
g^{(m)} \bigl( \omega, \omega_\alpha, (\omega^{j,i,1})_{j\leq p-1}\bigr)
$$
does not depend on $\omega_\alpha$ and the validity of the M\'ezard-Parisi ansatz follows as in the proof of the case $h=0$.  

If $q_*=0$ then we proved in Section \ref{Sec7label} that $\eta_i^\alpha = 0$ almost surely, so (\ref{Sec8eq1}) implies
\begin{equation}
\frac{1}{c-c^{-1}}
=
\frac{\e_{\alpha, i} t_1 s_{(1)} (1+ c t_1 s_{(1)})^{\zeta-1}}{\e_{\alpha,i} (1+ c t_1 s_{(1)})^\zeta}
\label{Sec8eq2}
\end{equation}
almost surely. Now, notice that when $q_*=0$, by Lemma \ref{Sec4Lem1} and (\ref{Sec6qstar}),
\begin{equation}
\e_{\alpha,i} s_{(1)} = 
\prod_{j\leq p-1} f^{(1)}(\omega, \omega_\alpha, \omega^{j,i,1})
= 
\prod_{j\leq p-1} f^{(1)}(\omega, \omega^{j,i,1}) =0
\end{equation}
almost surely. Therefore, the second derivative of the right hand side of (\ref{Sec8eq2}) at $t_1=0$ will be equal to $(\zeta-1)c\e_{\alpha,i} s_{(1)}^2$, since other terms will be equal to zero. Since the derivative of the left hand side of (\ref{Sec8eq2}) is zero, we get
\begin{equation}
\e_{\alpha,i} s_{(1)}^2 = 
\prod_{j\leq p-1} f^{(2)}(\omega, \omega_\alpha, \omega^{j,i,1}) = 0
\end{equation}
almost surely, which contradicts (\ref{Sec7Raa}). Therefore, $q_*$ can not be equal to zero in the presence of the external field, which finishes the proof of Theorem \ref{Th1}.
\qed

\end{document}